\documentclass{amsart}
\usepackage{amssymb}
\usepackage[matrix,arrow,curve]{xy}
\usepackage{verbatim}
\usepackage{a4wide,url}
\usepackage[applemac]{inputenc}
\usepackage{enumitem}

\def\Z{\mathbf{Z}}

\def\F{\mathbf{F}}
\def\Q{\mathbf{Q}}

\date{May 23, 2016.}
\usepackage{tikz}
\usetikzlibrary{backgrounds,fit}
\tikzset{every picture/.style = {inner sep=0pt,baseline}}
\tikzset{p/.style = {draw, shape          = circle,
                                 text           = black,
                                 font=\tiny,
                                 inner sep      = 0pt,
                                 outer sep      = 0pt,
                                 minimum size   = 1pt}}

\tikzset{r/.style = {red, very thin}}
\tikzset{g/.style = {green, very thin}}
\tikzset{b/.style = {blue, very thin}}
\tikzset{o/.style = {orange, very thin}}
\tikzset{u/.style = {black, very thin}}

\begin{document}
\centerline{}

\title{Homology computations for complex braid groups II}
\author[I.~Marin]{Ivan Marin}
\address{LAMFA, Universit\'e de Picardie-Jules Verne, Amiens, France}
\email{ivan.marin@u-picardie.fr}
\medskip

\begin{abstract} We complete the computation of the integral homology of the generalized
braid group $B$ associated to an arbitrary irreducible complex reflection group $W$ of
exceptional type. In order to do this we explicitely computed the recursively-defined differential
of a resolution of $\Z$ as a $\Z B$-module, using parallel computing. We also deduce from this general computation
the rational homology of the Milnor fiber of the singularity attached to most of these
reflection groups.
\end{abstract}

\maketitle

\section{Introduction}

This paper is a sequel of \cite{CALMAR}, where we managed to compute the
homology of many of the complex braid groups arising from complex reflection groups.
Parts of our difficulties in completing the project were computational in nature. Since then,
by using powerful computers during several months as well as a few computational tricks
we managed
 to complete some of the tables.
This paper is a report on these computations.

Recall from e.g. \cite{BMR} that one can attach to every finite complex reflection group $W$
a generalized braid group $B$. Without loss of generality, one can assume
that $W$ is irreducible. The usual
classification, due to Shephard-Todd, of
irreducible finite complex reflection groups, divides them into a general series
depending on 3 parameters, and a finite set of 34 exceptions. It is therefore
a natural question to explicitely compute the group homology of $B$ when $W$ belongs to
this finite set of exceptional groups ; moreover, one knows that $B$ has finite homological
dimension by work of Bessis in \cite{BESSISKPI1}.

By several arguments, recalled in \cite{CALMAR}, one can reduce the problem to
a fewer number of groups. In particular, the homology of groups of small rank
can be easily computed. Moreover, when $W$ is a \emph{real} reflection group,
a general complex due to Salvetti (see \cite{SALVETTI}) can be used to compute the homology of $B$
(see \cite{SALVETTI,CALMAR}); in this case, $B$ is an Artin group.

For these reasons, the remaining groups on which we need to focus are
the ones labelled $G_{24}$, $G_{27}$, $G_{29}$, $G_{31}$, $G_{33}$, $G_{34}$
in Shephard and Todd notation.
 Complexes can be obtained from the so-called
Garside theory introduced by Dehornoy and Paris. Indeed, Dehornoy and Lafont have
proven in \cite{DEHORNOYLAFONT} that a Charney-Meyer-Wittlesey-type complex can be used whenever $B$ is
the group of fraction of a so-called Garside monoid. Bessis has proven that,
when $W$ is well-generated, then $B$ satisfies this condition : there is one (and actually
several) convenient Garside monoid(s) $B^+$ (see \cite{BESSISKPI1}).
The one we use here has been specified for each group in \cite{CALMAR}.
 All the groups above are well-generated, \emph{but} $G_{31}$.
Even in the case of $G_{31}$, one can define a similar complex, by attaching to $B$ a
so-called Garside category instead of a Garside monoid. However, as we noticed in \cite{CALMAR},
the complexes obtained by this method are very big, which poses a computational memory
problem to compute their homology.

On the other hand, Dehornoy and Lafont introduced another, more mysterious but smaller kind of complex,
which can be attached to a similar Garside monoid, but for which a generalization to Garside categories
has not been proposed so far. Therefore, for this approach one needs to exclude the case of $G_{31}$.
The drawback of this complex is that, while the computation of the homology of the complex
is much easier as soon as it is explicitely described, the explicit computation of (the differential of) the complex itself is much more difficult
and time-consuming. For the other groups of rank at least 3, the specific Garside monoids chosen for these groups
have been specified in \cite{CALMAR}, table 7.

In the present work, we computed this differential. The result is stored inside large files, and could
in theory be used to compute the homology of $H_*(B,M)$ for an arbitrary $\Z B$-module $M$. In this
paper, we describe the result of $H_*(B,M)$ in the following cases :
\begin{enumerate}
\item $M = \Z$ with trivial action;
\item $M = \Z$ with action given by the sign morphism $B \to W \stackrel{\mathrm{det}}{\to} \{ \pm 1 \}$, which
exists because all the (pseudo-)reflections of $W$ have order $2$ in these cases;
\item $M = \Q[t,t^{-1}]$ with action given by the natural map $B \to  \Z$, $\sigma \mapsto t$.
\end{enumerate}
Note that, when $R$ is a commutative ring, $H_*(B,R[t,t^{-1}])$ can be identified with the homology of the Milnor fiber
of the singularity corresponding to $W$  (see \cite{CALARTIN}). For $G_{34}$, we were however unable to compute the homology
of the Milnor fiber in ranks $4$ to $6$ because of computer and software limitations.

It seems likely that the Dehornoy-Lafont complex
can be adapted to the kind of Garside categories
that are suitable for dealing with $G_{31}$, using its
description as a centralizer in the group of Coxeter type
$E_8$, as in \cite{BESSISKPI1}. However such a theory has not been
developped yet, and therefore $G_{31}$ is, for the time being, out of reach of our
computations.

\medskip

As an indication of computing time, we mention that the computation, for $G_{34}$, of the differential of 2000 of the 7414
5-cells lasted around 200 days of CPU time on a SMP architecture. The computation of the differentials
of the 5-cells and 6-cells altogether lasted around 78000 hours of CPU time.

\medskip

{\bf Acknowledgements.} The results presented here for $G_{34}$ were obtained using the ressources of the MeCS
platform of the Universit\'e de Picardie Jules Verne in Amiens, France. I thank very much Serge Van Criekingen for his help in using it.
I also thank Filippo Callegaro for several discussions and comments.

\section{Computational datas}

The size of the complexes we had to compute are tabulated in table \ref{tab:sizecomplexes}. Recall that each monoid $B^+$
is described as generated by a set $A$ of atoms, and that there is a distinguished element $\Delta$. Its set
$S$ of divisors is the same on the left and on the right, and is called the set of \emph{simples} of the monoid $B^+$. We describe the apparatus for the case $W$ is the complex reflection
group of type $G_{34}$, also called Mitchell's group.
Our programs and files are made to be primarily used by GAP3, but the syntax
is quite standard whence these files could be used by probably every
computer program with possibly only tiny changes to be made. These can be found
at \url{http://www.lamfa.u-picardie.fr/marin/G34homology.html}.

The group $W$ is described as a permutation group. The set $A$ is in 1-1 correspondence
with generators of $B^+$, which are stored in the file \verb+atomsG34.gap+ as an ordered list 
\verb+allatoms+ 
of 56 permutations. The set $S$ is in 1-1 correspondence with a subset of the
set of reflections of $W$, which is stored in the file  \verb+simplesG34.gap+
as an ordered list \verb+allsimples+ 
of 1584 permutations. An important additional data is the length of each simple as a product of atoms. This data is stored, in the same file, as the list \verb+simpleslengths+, in obvious bijection with
the list \verb+allsimples+, namely
\verb+simpleslengths[i]+ is the length ot the simple element \verb+allsimples[i]+.

By construction, the $k$-cells of the Dehornoy-Lafont complex
are in 1-1 correspondence with lists of atoms of the form $[a_1,\dots,a_k]$ with $a_i \in A$.
The files \verb+cells2N.gap+,\dots,\verb+cells6N.gap+ store them, under the variable name
\verb+cells2N+,\dots,\verb+cells6N+, as a list
of lists $[c_1,\dots,c_k]$ so that $a_i$ is the $c_i$-th atoms, namely $c_i$ is the position
in the list \verb+allatoms+ of the atom $a_i$. The 1-cells are simply given by the 56 atoms,
and the differential is simply $\partial( [a]) = a[\emptyset] - [\emptyset]$, where $[\emptyset]$
denotes the only 0-cell.

The files \verb+Dcells2.gap+,\dots,\verb+Dcells6.gap+ contain, under the variable names
\verb+Dcells2P+,\dots,\verb+Dcells6P+, the differentials of
the $k$-cells, for $k \in \{2,\dots, 6 \}$. The format is as follows. For example, the variable
\verb+Dcells4P+ is a list of 7520 elements $[v_1,\dots,v_{7520}]$,
where $v_r$ represents the differencial of the $r$-th cell in the list \verb+cells4+.
This differential is a linear combination of 3-cells with coefficients in the monoid algebra $\Z B^+$.
This linear combination can be written $\sum_{i=1}^{n} a_i b_i c_i$, with $a_i \in \Z$,
$b_i \in B^+$ and $c_i$ a 3-cell.The element $v_r$ is a list of $n$ terms
$[a_i,Q_i]$, where $Q_i= [\beta_i,\gamma_i]$
and $\beta_i, \gamma_i$ encode $b_i,c_i$. The encoding is as follows. An element of $B^+$
is a product of simple elements $s_1,\dots,s_q$. Then, $\beta_i$ is a list $[\sigma_1,\dots,\sigma_q]$ where $\sigma_j$ is the position of $s_j$ in the liste \verb+allsimples+. Finally,
$\gamma_i$ encodes the cell $c_i$ as in the file \verb+cells3N.gap+ (that is, as a list
of positions in the atom's list \verb+allatoms+).

\begin{table}
$$
\begin{array}{l|lrrrrrrr|}
\hline
 & & \mbox{0-cells} & \mbox{1-cells} & \mbox{2-cells} & \mbox{3-cells} & \mbox{4-cells} & \mbox{5-cells} & \mbox{6-cells} \\
\hline
G_{24} & CMW & 1 & 29 & 77 & 49 & & & \\
   & DL & 1 & 14 & 38 & 25  & & & \\
\hline
G_{27} & CMW & 1 & 41 & 115 & 75  & & & \\
   & DL & 1 & 20 & 62 & 43  & &  & \\
\hline
G_{29} & CMW & 1 & 111 & 635 & 1025 & 500 & & \\
   & DL & 1 & 25 & 127 & 207 & 104 & & \\
\hline
G_{33} & CMW & 1 & 307 & 3249 & 9747 & 11178 & 4374 &  \\
   & DL & 1 & 30 & 226 & 638 & 740 & 299& \\
\hline
G_{34} & CMW & 1 & 1583 & 31717 & 163219 & 337169 & 304927 & 100842 \\
   & DL & 1 & 56 & 711 & 3448 & 7520 & 7414 & 2686  \\
\hline
\end{array}
$$
\caption{Compared size of the complexes}
\label{tab:sizecomplexes}
\end{table}

\section{The results}

\begin{table}
$$
\begin{array}{|c|ccccccc|}
 & H_0 & H_1 & H_2 & H_3 & H_4 & H_5 & H_6 \\
\hline
G_{12} & \Z & \Z & 0 & & & & \\ 
G_{13} & \Z & \Z^2 & \Z & & & & \\ 
G_{22} & \Z & \Z & 0 & & & & \\
G_{24} & \Z & \Z & \Z & \Z & & & \\
G_{27} & \Z & \Z & \Z_3 \times \Z & \Z & & & \\
G_{29} & \Z & \Z & \Z_2 \times \Z_4 & \Z_2 \times \Z & \Z & & \\
G_{31} & \Z & \Z & \Z_6 & \Z & \Z & & \\
G_{33} & \Z & \Z & \Z_6 & \Z_6 & \Z & \Z & \\
G_{34} & \Z & \Z & \Z_6 & \Z_6 & \Z_3\times\Z_3\times\Z_6  & \Z_3\times\Z_3\times\Z & \Z \\
\end{array}
$$
\caption{Homology with trivial integer coefficients}\label{tabletriv}
\end{table}

\begin{table}
$$
\begin{array}{|c|ccccccc|}
 & H_0 & H_1 & H_2 & H_3 & H_4 & H_5 & H_6 \\
\hline
G_{12} &  \Z_2  &  \Z_3  & 0 & & & & \\ 
G_{13} &  \Z_2  &  \Z_2  &  \Z  & & & & \\ 
G_{22} &  \Z_2  &  0  & 0 & & & & \\
G_{24} &  \Z_2  &  0  &  \Z_2  &  0  & & & \\
G_{27} &  \Z_2  &  0  &  \Z_2    &  0  & & & \\
G_{29} &  \Z_2  &  0  &  \Z_2 \times \Z_{4}    &  \Z_2 \times \Z_{40}    &  0  & & \\
G_{33} &  \Z_2  &  0  &  \Z_2  &  \Z_2  &  \Z_2  &  0  & \\
G_{34} &  \Z_2  &  0  &  \Z_6  & \Z_2 & \Z_6 & \Z_{252} & 0 \\
\end{array}
$$
\caption{Homology with coefficients in the integer sign representation}\label{tablesign}
\end{table}

\subsection{Hand-made computations for groups of rank 2}

For the groups $G_{12}$ and $G_{22}$ we use the monoids
$$
\begin{array}{lcl}
B_{12}^+ &=& \langle x_1,x_2,x_3 \ \mid \ x_1 x_2 x_3 x_1 = x_2 x_3 x_1 x_2 = x_3x_1 x_2 x_3 \rangle \\
B_{22}^+ &=& \langle x_1,x_2,x_3 \ \mid \ x_1 x_2 x_3 x_1x_2 = x_2 x_3 x_1 x_2x_3 = x_3x_1 x_2 x_3x_1 \rangle \\
\end{array}
$$
It is proved in \cite{PICANTIN} (see ch. 5 ex. 11 and ch. 6) that these monoids are Garside
monoids, so we can build a Dehornoy-Lafont complex for each one of them.
In both cases, the set of atoms is $\mathcal{X} = \{ x_1,x_2,x_3 \}$, which we endow with
the linear ordering $x_1 < x_2 < x_3$. The 2-cells are then $[x_1,x_2]$ and $[x_1,x_3]$
and there are no cells of higher degree. Therefore the only datas to be computed in order to
define the corresponding chain complex are
$\partial_2 [x_1,x_2]$ and $\partial_2 [x_1,x_3]$. For $G_{12}$ one gets
$$
\begin{array}{lcl}
\partial_2 [x_1,x_2] &=& x_2 x_3 x_1 [x_2] - x_1 x_2 x_3 [x_1] - x_1 x_2 [x_3] - x_1[x_2]
- [x_1] + x_2 x_3 [x_1] + x_2 [x_3] + [x_2] \\
\partial_2 [x_1,x_3] &=& x_3 x_1 x_2[x_3] - x_1 x_2 x_3 [x_1] - x_1 x_2 [x_3] - x_1[x_2] - [x_1]
+x_3 x_1 [x_2] + x_3[x_1] + [x_3]
\end{array}
$$
and for $G_{22}$ one gets
$$
\begin{array}{lcr}
\partial_2 [x_1,x_2] &=& x_1 x_2 x_3 x_1[x_2]-x_3x_1x_2x_3[x_1] - x_3 x_1 x_2[x_3] - x_3 x_1[x_2]
-x_3[x_1]-[x_3] + x_1x_2x_3[x_1] \\ & &+ x_1x_2[x_3] + x_1[x_2] + [x_1] \\
\partial_2 [x_1,x_3] &=& x_2 x_3x_1x_2[x_3] - x_3 x_1 x_2 x_3[x_1]-x_3x_1x_2[x_3] - x_3x_1[x_2]
-x_3[x_1]-[x_3]+x_2x_3x_1[x_2] \\ & &+ x_2 x_3[x_1]+x_2[x_3]+[x_2]
\end{array}
$$
The $\Z B$-module structure on $\Q[t,t^{-1}]$ affording the homology of the Milnor fiber
is given by $x_i \mapsto t$ in both cases (see \cite{CALCLASS, CALARTIN} for the classical case).

For the group $G_{13}$, we use the fact that its braid group is isomorphic to the Artin group of type
$I_2(6)$. More precisely, in this case $B$ has for presentation $\langle x,y,z \ | \ yzxy=zxyz,
zxyzx=xyzxy \rangle$ and an isomorphism with the Artin group $\langle a,b \ | \ ababab=bababa \rangle$
is given by the formulas
$$
\left\lbrace \begin{array}{lcl} 
a &=& zx \\
b &=& zxy(zx)^{-1}
\end{array}
\right.
\ \ \ \left\lbrace \begin{array}{lcl} 
x &=& (baba)^{-1} abaa \\
y &=& a^{-1}ba \\
z &=&  (aba)^{-1}b(aba)
\end{array}
\right.
$$
In particular we have $B^{ab} \simeq \Z^2$, and two different bases
are afforded by $(\bar{x}, \bar{y} = \bar{z})$ and $(\bar{a},\bar{b})$,
the relation between them being $\bar{x} = \bar{a}- \bar{b}$, $\bar{y} = \bar{b}$.
The (homological) Salvetti complex (see \cite{CONCSALV}) is made of free $\Z B$-modules
with basis $[\emptyset]$ (one 0-cell), $[a],[b]$ (two 1-cells) and $[a,b]$ (one 2-cells).
The differential is given by $\partial [\emptyset] = 0$,
$\partial [a] = (a-1)[\emptyset]$, $\partial [b] = (b-1)[\emptyset]$, and
$$
\partial[a,b] = (1-a+ab-aba+abab-ababa)[b] - (1-b+ba-bab+baba-babab)[a].
$$
Therefore, tensoring by the $\Z B$-module $R=\Z[t,t^{-1}]$-module defined
by $x,y,z \mapsto t$, we get the differential
$d[a] = (t^2-1)[\emptyset]$,
$d[b] = (t-1)[\emptyset]$,
$d[a,b] = (1+t^3+t^6)(1-t) \left( (1+t)[b] - [a]\right)$.
Specializing at $t=1$ and $t = -1$ we readily get the results of tables \ref{tabletriv} and \ref{tablesign} for $G_{13}$.
In general, one gets the following homology of the Milnor fiber for $G_{13}$ :
$$
H_0 = R/(t-1)R \simeq \Z, 
H_1 = R/(1-t)(1+t^3+t^6) \simeq \Z^6,
H_2 = 0.$$

\subsection{Milnor fibers}

We let $\Phi_n \in \Z[t]$ denote the $n$-th (rational) cyclotomic polynomial, and $R = \Q[t,t^{-1}]$. In table 
\ref{tablemilnor} we indicate the homology $H_*(B,\Q[t,t^{-1}]) $ as a $R$-module, which can be identified with the rational homology of the
Milnor fiber. Results on the classical (Artin) cases can be found in \cite{CALARTIN,CALCLASS}. In our paper the $R B$-module structure
on $\Q[t,t^{-1}]$ is given by $\sigma_i \mapsto t$ (while the choice in \cite{CALCLASS} is $\sigma_i \mapsto -q$).
This homology was computed using the software \verb+Macaulay2+, see \cite{MACAULAY2}.
In the table, for each $P \in R$, the presence of  $P$ in the table symbolizes the $R$-module $R/(P)$, and $\Q$ is a shortcut for $R/\Phi_1 = \Q[t,t^{-1}]/(t-1)$,
that we use for $H_0$ and $H_1$.
Notice that
$$
\frac{t^{20}-1}{t+1} = \frac{t^{20}-1}{\Phi_2}=\Phi_1\oplus\Phi_4\oplus\Phi_5\oplus\Phi_{10}\oplus\Phi_{20}  
$$
It follows from the table that the Poincaré polynomials of the Milnor fibers are as follows.
$$
\begin{array}{|l|l|}
\mbox{Group} & \mbox{Poincar\'e polynomial} \\
\hline
G_{12} & 1+x+6x^2\\
G_{13} & 1+6x\\
G_{22} & 1+x+8x^2\\
G_{24} & 1+x+9x^3 \\
G_{27} & 1+x+17x^3 \\
G_{29} & 1+x+4x^3+21x^4 \\
G_{33} & 1+x+13x^5 \\
G_{34} & 1+x+2x^3+\dots ?\\
\hline
\end{array}
$$

\begin{table}
$$
\begin{array}{|c|ccccccc|}
 & H_0 & H_1 & H_2 & H_3 & H_4 & H_5 & H_6 \\
\hline
G_{12} &  \Q  &  \Q  & \Phi_{6}\oplus\Phi_{12} & & & & \\ 
G_{13} &  \Q  &  \Q\oplus \Phi_9  &  0  & & & & \\ 
G_{22} &  \Q  &  \Q  & \Phi_{15} & & & & \\
G_{24} &  \Q  &  \Q  &  0  &  \Phi_1\oplus\Phi_3\oplus \Phi_7  & & & \\
G_{27} &  \Q  &  \Q  &  0    &  (t^{15}-1)\oplus \Phi_3  & & & \\
G_{29} &  \Q  &  \Q  &  0    &  \Phi_4\oplus\Phi_4    &  \frac{t^{20}-1}{t+1}\oplus \Phi_4 & & \\
G_{33} &  \Q  &  \Q  &  0  &  0  &  0  &  (t^9-1) \oplus \Phi_5  & \\
G_{34} &  \Q  &  \Q  &  0  & \Phi_6 & ? & ? & ? \\
\end{array}
$$
\caption{Rational homology $H_*(B,\Q[t,t^{-1}])$ of the Milnor fiber}\label{tablemilnor}
\end{table}

With the same algorithm, for any given $p$ we can compute the homology $H_*(B,\F_p[t,t^{-1}])$ modulo $p$
of the Milnor fiber. We compute this for $p \in \{2,3,5,7 \}$. These numbers might be particularly interesting
for applications because they are the only primes dividing $|W|$ for $W$ in our list. Letting $\Phi_n$ denote the reduction
modulo $p$ of the $n$-th usual (rational) cyclotomic polynomial, we get the same result as in table \ref{tablemilnor},
\emph{except} for the following cases :
\begin{itemize}
\item When $W = G_{12}$ and $p \in \{2,3,5,7\}$ we have $H_2(B,\F_p[t,t^{-1}]) = P_{12,2}$ with $P_{12,2} = t^6-t^5+t^3-t+1$.
When $p =3,5,7$ this amounts to $\Phi_6 \oplus \Phi_{12}$ as in the rational case, but in case $p=2$ we have $P_{12,2} = (t^2+t+1)^3 = \Phi_3^3$.
\item When $W = G_{13}$ and $p \in \{ 2,3,5,7 \}$ we have $H_1(B,\F_p[t,t^{-1}]) = P_{13,1}$ with $P_{13,1} = (1-t)(1+t^3+t^6)$.
When $p =2,5,7$ this amounts to $\Phi_1 \oplus \Phi_{9}$ as in the rational case, but in case $p=3$ we have $P_{13,1} = (t-1)^7$.
\item When $W = G_{24}$, and $p \in \{2,3,5,7 \}$, we have $H_3(B,\F_p[t,t^{-1}]) = P_{24,3}$
with $P_{24,3}(t) = t^9+t^8+t^7-t^2-t-1$. When $p=2,5$ this amounts to $\Phi_1\oplus \Phi_3 \oplus \Phi_7$
as in the rational case, but in case $p = 3$ we have $P_{24,3} = (t-1)^3 \Phi_7$
while for $p=7$ we have $P_{24,3} = (t+3)(t-1)^7(t+5)$.
\item When $W = G_{29}$, and $p=2$, we have $H_3(B,\F_p[t,t^{-1}]) = (t+1)^3 \oplus \Phi_4$ (instead of $\Phi_4 \oplus \Phi_4$).
\item When $W = G_{33}$, and $p=2,3$, we have $H_3(B,\F_p[t,t^{-1}]) =\Phi_1$ (instead of $0$).
\item When $W = G_{34}$, and $p=2,3$, we have $H_3(B,\F_p[t,t^{-1}]) =\Phi_6 \oplus \Phi_1$ (instead of $\Phi_6$).
\item When $W = G_{29}$, and $p=2$, we have $H_4(B,\F_2[t,t^{-1}]) =(t^{20}-1)\oplus \Phi_4 $ (instead of $\frac{t^{20}-1}{t+1} \oplus \Phi_4$).
\item When $W = G_{33}$, and $p=2,3$, we have $H_4(B,\F_p[t,t^{-1}]) =\Phi_1 $ (instead of $0$).
\item When $W = G_{33}$, and $p \in \{2,3,5,7 \}$, we have $H_5(B,\F_p[t,t^{-1}]) =P_{33,5}$
with $P_{33,5} = t^{13}+t^{12}+t^{11} + t^{10} + t^9 - t^4-t^3-t^2-t-1$. When $p =2,3,7$ this amounts to
$(t^9-1)\oplus \Phi_5$ as in the rational case, but when $p = 5$ we have $P_{33,5} = (t-1)^5(t^2+t+1)(t^6+t^3+1)$.
\end{itemize}

F. Callegaro observed and communicated to us that some of the special cases we obtain here can be unified by using that
$\Phi_{mp^i}(t) = \Phi_m(t)^{\phi(p^i)} \mod p$ (see e.g. \cite{GUERRIER}). Moreover,
it follows from these computations that there are no torsion of order $2,3,5,7$ in the integer
homology of the Milnor fiber for  $G_{12}, G_{13}, G_{22}, G_{24}, G_{27}$, while there is 2-torsion for $G_{29}$
and 2,3-torsion for $G_{33}$. 

We were not able to compute the higher dimensional rational homology groups of the Milnor fiber for $G_{34}$, because
of memory issues with the software we are using. However, someone having more computational skills than
us might well be able to compute them, using smarter procedures. Therefore, we included the \verb+Macaulay2+ script
among the datafiles which can be found on our webpage.

\end{document}